\documentclass[a4paper,11pt]{amsart}
\usepackage{amssymb}

\DeclareMathSymbol{\subsetneqq}{\mathbin}{AMSb}{36}

\DeclareMathSymbol{\subsetneqq}{\mathbin}{AMSb}{36}

\newcommand{\R}{\mathbb{R}}

\newcommand{\dint}{\displaystyle\int}

\newtheorem{thm}{{\bf Theorem}}[section]
\newtheorem{lem}{{\bf Lemma}}[section]
\newtheorem{prop}{{\bf Proposition}}[section]

\theoremstyle{remark}
\newtheorem{rem}{\bf Remark}[section]

\theoremstyle{definition}
\newtheorem{defi}{\bf Definition}[section]

\author{Olfa Bejaoui}
\address{Faculty of Sciences of Tunis,
Department of Mathematics, ElManar 2092, Tunis, Tunisia.}
\email{ bjaouiolfa@yahoo.fr}

\title[Axisymmetric Rotating Fluid Equations]
{Axisymmetric Rotating Fluid Equations}
\date{\today}

\begin{document}
\begin{abstract}
We investigate the equations of anisotropic axisymmetric incompressible viscous
fluids in the exterior of a cylinder of $\R^3$, rotating around an inhomogeneous vector
  $B( t, r)$. We prove uniform local existence with respect to the Rossby number in suitable anisotropic Sobolev spaces. We also obtain the propagation of the isotropic Sobolev regularity. This extends the results of \cite{paicu'}.
\end{abstract}
%@@@@@@@@@@@@@@@@@@@@@@@@@@@@@@@@@@@@%@@@@@@@@@@@@@@@@@@@@@@@@@@@@@@@@@@@@%@@@@@@@@@@@@@@


%@@@@@@@@@@@@@@@@@@@@@@@@@@@@@@@@@@@@%@@@@@@@@@@@@@@@@@@@@@@@@@@@@@@@@@@@@%@@@@@@@@@@@@@@

\maketitle

%@@@@@@@@@@@@@@@@@@@@@@@@@@@@@@@@@@@@%@@@@@@@@@@@@@@@@@@@@@@@@@@@@@@@@@@@@%@@@@@@@@@@@@@@@
\section{Introduction}
The motion of incompressible rotating fluids in a domain $\Omega$ of
$\mathbb{R}^3$ is described by the following system of equations
\begin{eqnarray*}(S^\varepsilon)
\begin{cases}
\partial_t u^\varepsilon+ (u^\varepsilon.\nabla) u^\varepsilon-\nu_h\Delta_h
u^\varepsilon-\nu_v \partial^2_{x_3} u^\varepsilon
+\dfrac{1}{\varepsilon}\, ( u^\varepsilon \times B)+\nabla
p^\varepsilon=0\,\,\,
\mbox{in}\,\,\,\Omega,\\
\textmd{div }{u}^\varepsilon=0\,\,\,\mbox{in}\,\,\, \Omega,\\
u^\varepsilon=0\,\,\, \mbox{on}\,\,\,\partial\Omega,
\\u^\varepsilon(0,x)=u_0(x),\end{cases}
\end{eqnarray*}
where $u^\varepsilon$ is the velocity field and $p^\varepsilon$ is
the pressure. The constants $\nu_h>0$ and $\nu_v\geq 0$ represent
respectively the horizontal and vertical viscosities. The symbol
$\Delta_h$ stands for the horizontal Laplacian and the term
$\dfrac{1}{\varepsilon}\,( u^\varepsilon \times B)$ represents the
Coriolis force, where $B$ is the rotation vector and $\varepsilon$
is a small parameter. We assume that $B$ is a smooth function with
bounded derivatives, depending on time $t$ and horizontal variables
$x_h$ (that is $x=(x_h, x_3)$). Generally, it is a vector field
directed along the $x_3$ coordinate. Additional assumptions on $B$ will be made later.\\

Notice that if $B=0$ in the system $(S^\varepsilon)$, then we get
the classical incompressible Navier-Stokes equations. It is well known (see \cite{leray}) that if $u_0$ is
only in $L^2(\mathbb{R}^3)$ and  $\nu_v>0$, then a global weak solution exists. The uniqueness of such solution is an outstanding open problem. Concerning strong solutions, the pioneer work goes back to Fujita-Kato \cite{FK} where local existence and uniqueness were obtained in $\dot{H}^{\frac{1}{2}}(\mathbb{R}^3)$. Moreover, if the initial data is small enough then the solution is global in time. We also refer to \cite{Iftimie} for well-posedness in thin domains and in the anisotropic Sobolev spaces $H^{0,s}$. We recall that $H^{0,s}$ is the space of functions which
are $L^2$ in the horizontal variables and $H^s$ in the vertical one.\\

Let us then consider the case $\nu_v=0$. The anisotropic
Navier-Stokes system with vanishing vertical viscosity was studied
for the first time in \cite{CDGG}, where local existence for large data and global existence for small data were obtained in the anisotropic
Sobolev spaces $H^{0,s}(\mathbb{R}^3)$, $s>\frac{1}{2}$. Note that in \cite{CDGG} the uniqueness was proved only for $s>\frac{3}{2}$. Later on, D. Iftimie (\cite{Iftimie2}) filled the gap between existence and uniqueness by proving uniqueness for $s>\frac{1}{2}$. Recently, M.
Paicu obtains the uniqueness in the critical Besov space
${\mathcal B}^{0,\frac{1}{2}}(\mathbb{R}^3)$ {\footnote{This space is close to $H^{0,\frac{1}{2}}$. However, as far as we know, there is no result in $H^{0,\frac{1}{2}}$.}}(see \cite{paicu}).\\

Let us now
recall some well-known facts about the constant case $B=e_3:=(0,0,1)$. For a physical motivation we refer the reader to the book of J. Pedlosky
\cite{pedlosky} as well as to \cite{CDGG-Book, GSR-HB}. As the singular perturbation
is a linear skew-symmetric operator, weak solutions can be constructed by the approximate scheme of Friedrichs when  $\nu_v>0$:
approximate solutions are obtained by a standard truncation in high frequencies.
In \cite{CDGG}, J.-Y. Chemin {\it{et al.}} obtained local existence in
the anisotropic Sobolev spaces $H^{0,1/2+\epsilon}(\R^3)$ and the global
existence for data which are small compared to the horizontal
 viscosity.  They also proved the global existence of the solution for
anisotropic rotating fluids. In \cite{paicu1},
global existence of the solution for rotating fluids with vanishing
vertical viscosity was shown in the periodic case. Other related results can be found in \cite{CDGG, CDGG1, paicu1}.\\

Let us focus now on the variable case $B=B(t,x_h)$. As for the constant case, the classical proofs of existence of weak
solutions for the Navier-Stokes equations can be
extended to $(S^\varepsilon)$ when $\nu_v>0$. The asymptotic of those solutions was investigated by I. Gallagher and L. Saint-Raymond in \cite{GR}. Using weak compactness arguments, they showed that weak solutions
converge to the solution of a heat equation in the region when $B$ is non stationary.\\

The existence of
strong solutions in Sobolev spaces was the main goal of a recent work of M.
Majdoub and M. Paicu \cite{paicu}. They
obtained global existence for small initial data, and uniform local existence for
large data under the assumption that the field $B$ depends on $t$
and $ x_1$ (or $t$ and $x_2$).\\

 In this paper, our main concern is to improve this assumption on the field $B$ in order to get uniform local  existence in the general case. To do so, we restrict
ourselves to the axisymmetric case. This means that we assume that the velocity $u^{\varepsilon}$ and the pressure $p^{\varepsilon}$ are axisymmetric (see Definition 2.3 below). We also assume that the domain $\Omega$ is the exterior of some cylinder. We obtain uniform local existence with respect to the Rossby number ${\varepsilon}$ as well as the propagation of the isotropic Sobolev regularity.\\
\indent It is
expected that similar results can be shown in the case of
the whole space $\R^3$. This will be dealt with in a forthcoming work.\\

The paper unfolds as follows: section $2$
contains some notations needed in the statement and the proofs of our
results. In the third section, we present the functional spaces used along this paper and we state the
main results. Section $4$ is devoted to the proof of the uniform local existence and to the propagation of  the isotropic Sobolev regularity.
A few technical lemmas have been postponed in the final section.\\ \indent Finally, $C$ will denote a constant that does not depend on
$\varepsilon$ but that may change from line to line.

%@@@@@@@@@@@@@@@@@@@@@@@@@@@@@@@@@@@@%@@@@@@@@@@@@@@@@@@@@@@@@@@@@@@@@@@@@%@@@@@@@@@@@@@@@
\section{Definitions and notations}
Here we give the definitions of axisymmetric domains and
axisymmetric vector fields.
\begin{defi}
We denote by $(r,\theta,x_3)\in \mathbb{R}^*_+\times
]-\pi,\pi[\times \mathbb{R}$ the cylindrical coordinates in
$\mathbb{R}^3$, where $r$ and $\theta$ are defined
by$$r=\sqrt{x_1^2+x_2^2},\,\,x_1=r\cos \theta \,\,\mbox{and}\,\,
x_2=r \sin\theta.$$\end{defi}
\begin{defi}
An open domain $\Omega$ of $\mathbb{R}^3$ is said to be axisymmetric
if for every rotation $\Gamma$ around the vertical axis $e_3$, we
have $\Gamma(\Omega)\subset \Omega$.\end{defi}
\begin{defi}
An axisymmetric vector field $u$, defined on an open axisymmetric
domain of $\mathbb{R}^3$, is a field having the
representation{\footnote{ $u^\theta$ is called the swirl component. If $u^\theta=0$, we say that $u$ is an axisymmetric vector field without swirl.}}$$u(r,x_3)=u^r(r,x_3)e_r+u^\theta(r,x_3)e_\theta+u^3(r,x_3)e_3$$in
the cylindrical coordinate system,
where$$e_r=(\frac{x_1}{r},\frac{x_2}{r},0),\,\,e_\theta=(-\frac{x_2}{r},\frac{x_1}{r},0),\,\,e_3=(0,0,1
).$$\end{defi}

Let us then introduce some useful notations. We denote by
$$\mathcal{V}(\Omega)=\Big\{\;u\in
\left(C^\infty_c(\Omega)\right)^3,\quad\textmd{div}u=0\;\,\Big\},$$ where $C^\infty_c(\Omega)$
is the space of smooth functions compactly supported in the domain
$\Omega:=\Omega_h\times
\mathbb{R}$. Here $\Omega_h$ is given by
$$
\Omega_h:=\left\{\;(x_1,x_2)\in \mathbb{R}^2,\quad
\,x_1^2+x_2^2>\rho^2>0\;\right\},$$ where $\rho$ is a fixed positive number.
 The notation $(./.)_E$ corresponds to the inner product
in the functional space $E$ and the symbol $\partial_i$
stands for the partial derivative in the direction $x_i$.\\For
an axisymmetric vector field $u$ defined on $\Omega$,
we write $u(x_h,x_3)=\tilde{u}(r,x_3)$, and we define its vertical Fourier
transform by
$$\mathcal{F}^\vee(u)(x_h,\xi_3)
\,=\,\int_{\mathbb{R}} u(x_h,x_3)\, e^{-i \xi_3.x_3}\,dx_3.$$ The
operator of localization in vertical frequencies $S_N^{x_3}$ ($N\in\mathbb{N}$), is defined by
$$\mathcal{F}^{\vee}(S_N^{x_3}u)(x_h,.)=\Psi\left(2^{-N}|.|\right)\mathcal{F}^{\vee}(u)(x_h,.),$$
where $\Psi$ is a smooth compactly supported function with
values in [0,1] such that
\begin{eqnarray*}
\begin{cases}
\Psi(s)=1,\quad \mbox{if}\quad s\in [0,1]\\ \Psi(s)=0,\quad
\mbox{if}\quad |s|\geq 2.
\end{cases}
\end{eqnarray*}

\section{Functional spaces and statement of the results}
In order to proceed in a more easy way, we give definitions and properties of some functional spaces used
along this paper.
In the frame of anisotropic Lebesgue spaces, the H\"older inequality
reads.
\begin{lem}\label{HY}
{\bf{H\"older inequality}}\quad \\
Let $1\leq p,\,p',\, p'',\,q,\,q',\,q''\leq \infty$ be numbers such
that $\displaystyle{\frac{1}{p} = \frac{1}{p'} +
\frac{1}{p''}}\,\mbox{and}\,\displaystyle{\frac{1}{q} = \frac{1}{q'}
+ \frac{1}{q''}}$. Then, we have
$$\|u v\|_{L^p_v(L^q_h)} \leq \|u\|_{L^{p'}_v(L^{q'}_h)}\|v\|_{L^{p''}_v(L^{q''}_h)}.$$\end{lem}
Let us recall the definition of isotropic Sobolev
spaces.\begin{defi}Let $m\geq 1$ be an integer. The space
$\mathcal{H}^m_0(\Omega)$ is defined as the closure of
$\mathcal{V}(\Omega)$ in the space $H^m(\Omega)$ for the norm $\|.\|_{H^m(\Omega)}$.
\end{defi}

Considering the anisotropy of the problem, we use spaces of
functions that take into account this anisotropy. More precisely, we
use anisotropic Sobolev spaces. Such spaces have been introduced by
D. Iftimie in \cite{Iftimie}.
\begin{defi}
Let $s$ be a real number. The norm $\|.\|_{H^{0,s}(\Omega)}$ is
defined by
$$\|u\|_{H^{0,s}
(\Omega)}^2=\int_{\Omega_h}\|u(x_h,.)\|_{H^s_{\vee}}^2\,d x_h,$$
where
$$\|u(x_h,.)\|_{H^s_{\vee}}^2=\int_{\mathbb{R}}
(1 + \xi_{3}^2)^{s} |\mathcal{F}^{\vee}(u)(x_h,\xi_3)|^2\,d\xi_3.$$
We define the norm $\|.\|_{H^{1,s}(\Omega)}$ by
$\|u\|_{H^{1,s}(\Omega)}=\|\nabla_h u\|_{H^{0,s}(\Omega)}$.
\end{defi}
Throughout this paper, we consider spaces constructed on the Sobolev
space $H^{0,s}(\Omega)$. \begin{defi}The space $
\mathcal{H}_0^{0,s}(\Omega)$ is the closure of $\mathcal{V}(\Omega)$ for the $\|.\|_{H^{0,s}(\Omega)}$ norm.
The space $ \mathcal{H}_0^{1,s}(\Omega)$ is the closure of
$\mathcal{V}(\Omega)$ for the $\|.\|_{H^{1,s}(\Omega)}$
norm.\end{defi}

Let us stress that in all what follows, we consider $\Omega=\Omega_h\times
\mathbb{R}$ and we suppose that $B=B(t,r)$ is a $C^\infty$
vector field defined on $\Omega_h$.\\

Now, we are ready to state the
main results of this paper.
\begin{thm}\label{th5}Assume that $u_0\in \mathcal{H}_0^{0,s}(\Omega)$ is an
axisymmetric vector field with $s>\frac{1}{2}$.
Then, there exists a time $T>0$ independent of $\varepsilon$ and a unique solution $u^\varepsilon$ to $(S^\varepsilon)$ such that $$u^\varepsilon \in{\mathcal
C}([0,T], \mathcal{H}_0^{0,s}(\Omega))\cap L^2([0,T],
\mathcal{H}_0^{1,s}(\Omega)).$$
\end{thm}
\begin{thm}\label{th6}Let $m\geq 1$ be an integer and assume that $u_0\in \mathcal{H}^m_0(\Omega)$ is an
axisymmetric vector field. Then, there exists a time $T>0$ independent of $\varepsilon$ and a unique solution $u^\varepsilon$ to $(S^\varepsilon)$ such that $$u^\varepsilon \in{\mathcal
C}([0,T], \mathcal{H}^m_0(\Omega));\,\nabla_h u^\varepsilon \in
L^2([0,T], \mathcal{H}^m_0(\Omega)).$$
\end{thm}
We mention that Theorem \ref{th6} is a consequence of Theorem
\ref{th5} and of the following result about the Navier-Stokes
equations with vanishing vertical viscosity proved in
\cite{paicu'} in the case of the whole space.

\begin{thm}
\label{th1}Let
$m\geq 1$ be an integer and assume that $u_0$ belongs to $
\mathcal{H}^m_0(\Omega)$. Then, there exists a positive time $T$
such that the Navier-Stokes
equations with vanishing vertical viscosity admits a unique solution $u$ satisfying
$$u\in{\mathcal C}([0,T],\mathcal{H}^m_0(\Omega))\,;\,\nabla_h
u\in L^2([0,T],\mathcal{H}^m_0(\Omega)).$$
 Let us denote by $T^*$ the maximal time of
existence; if $T^*$ is finite, then
 $$\lim_{t \rightarrow
T^*}\int_0^t\|\nabla_h
u(\tau)\|_{L^\infty_v(L^2_h)}^2\big(1+\|u(\tau)\|^2_{L^\infty_v(L^2_h)}\big)d\tau=+\infty.$$
\end{thm}
%@@@@@@@@@@@@@@@@@@@@@@@@@@@@@@@@@@@@%@@@@@@@@@@@@@@@@@@@@@@@@@@@@@@@@@@@@%@@@@@@@@@@@@@@@
\section{Uniform local existence results}
The main goal of this section is to prove uniform local
existence of strong solutions with respect to the Rossby number $\varepsilon$. The first step consists in
splitting the initial data into a small part and a regular one.
Hence, we get two systems: a globally well-posed linear system
associated to the regular part of $u_0$ with solution $v^\varepsilon$ and a nonlinear one
associated to the small part with solution $w^\varepsilon:=u^\varepsilon-v^\varepsilon$. We have only to show the
local uniform existence of $w^\varepsilon$. As in \cite{paicu'}, we use energy
estimate and a Gronwall's lemma. The most delicate term to estimate is
$(v^\varepsilon.\nabla
v^\varepsilon/w^\varepsilon)_{H^{0,s}}$, and here comes the
importance of considering axisymmetric vector fields and the fact that the
domain $\Omega$ does not contain a neighborhood of zero. The
key idea is, then, to use extension operators and Sobolev
embeddings.\\

We now come to the details of the proof of Theorem \ref{th5}.

\begin{proof}[Proof of Theorem \ref{th5}]\quad\\
First of all, using the operator of localization in vertical
frequencies $S_N^{x_3}$, we decompose $u_0$ into two parts. As $u_0$
belongs to the space $\mathcal{H}_0^{0,s}(\Omega)$, we
obtain$$\lim_{N'\rightarrow +\infty}
\|(I-S_{N'}^{x_3})u_0\|_{H^{0,s}}=0.$$Hence, there exists a positive
integer $N$ such that $\|(I-S_{N}^{x_3})u_0\|_{H^{0,s}}\leq c
\nu_h$, where $c>0$ is a small constant. Then, we split the system
$(S^\varepsilon)$ into
\begin{eqnarray*}(S_1^\varepsilon)
\begin{cases}
\partial_t v_N^\varepsilon-\nu_h\Delta_h
v_N^\varepsilon+\dfrac{1}{\varepsilon}\, ( v_N^\varepsilon\times
B)+\nabla p_N^\varepsilon=0\,\,\,
\mbox{in}\,\,\,\Omega,\\
\textmd{div}\,v_N^\varepsilon=0\,\,\,\mbox{in}\,\,\,\Omega,\\
v^\varepsilon_{N}=0\,\,\,\mbox{}on\,\,\,\partial\Omega,\\v^\varepsilon_N(0,x)=S_N^{x_3}
u_0,
\end{cases}
\end{eqnarray*}
and
\begin{eqnarray*}(S_2^\varepsilon)
\begin{cases}
\partial_t w_N^\varepsilon+ (w_N^\varepsilon+v_N^\varepsilon).\nabla (w_N^\varepsilon+v_N^\varepsilon)-\nu_h\Delta_h
w_N^\varepsilon+\dfrac{1}{\varepsilon}\, (w_N^\varepsilon\times
B)+\nabla p_N^\varepsilon=0\,\,\,
\mbox{in}\,\,\,\Omega,\\
\textmd{div}\,w_N^\varepsilon=0\,\,\,\mbox{in}\,\,\,\Omega,\\
w^\varepsilon_{N}=0\,\,\,\mbox{on}\,\,\,\partial\Omega,\\w^\varepsilon_{N}(0,x)=(I-S_N^{x_3})u_0.
\end{cases}
\end{eqnarray*}
Let us notice that $S_N^{x_3}u_0$ belongs to
$\mathcal{H}_0^{0,s}(\Omega)$. As $(S_1^\varepsilon)$ is a linear
system with a regular initial data, there exists a unique global in
time solution$$v_N^\varepsilon \in L^\infty(\mathbb{R}_+,
\mathcal{H}_0^{0,s}(\Omega)) \cap L^2(\mathbb{R}_+,
\mathcal{H}_0^{1,s}(\Omega)).$$Therefore, we need only to prove the
uniform local existence of $w_N^\varepsilon$. In other words, we
will prove that $w_N^\varepsilon$ which is small with respect to
$\nu_h$ at $t=0$  remains so for a certain time. Thus, let us
define
$$T_{\varepsilon,N}=\sup\Big\{t \geq
0/\,\,\forall \,\,0\leq t'\leq t,\quad
\|w_N^\varepsilon(t')\|_{H^{0,s}}\leq 2 c \nu_h\Big\}.$$As
$w_N^\varepsilon$ is a divergence free vector field, we
obtain$$(\nabla p_N^\varepsilon/w_N^\varepsilon)_{H^{0,s}}=-(
p_N^\varepsilon/\, \textmd{div}\,w_N^\varepsilon)_{H^{0,s}}=0,$$ and
as the vector field $B$ is independent of $x_3$, we get $$(
w_N^\varepsilon\times B/w_N^\varepsilon)_{H^{0,s}}=0.$$Hence,
computing the $H^{0,s}$ scalar product of $(S_2^\varepsilon)$ by
$w_N^\varepsilon$ leads to
$$\frac{1}{2}\frac{d}{dt}\|w_N^\varepsilon\|_{H^{0,s}}^2+ \nu_h
\|\nabla _h w_N^\varepsilon\|_{H^{0,s}}^2\leq
\widetilde{T}^{\varepsilon,\,1}_N+\widetilde{T}^{\varepsilon,\,2}_N+\widetilde{T}^{\varepsilon,\,3}_N+
\widetilde{T}^{\varepsilon,\,4}_N,$$where\begin{eqnarray*}
\widetilde{T}^{\varepsilon,\,1}_N&=&|(w_N^\varepsilon.\nabla
w_N^\varepsilon/w_N^\varepsilon)_{H^{0,s}}|,\\
\widetilde{T}^{\varepsilon,\,2}_N&=&|(v_N^\varepsilon.\nabla
w_N^\varepsilon/w_N^\varepsilon)_{H^{0,s}}|,\\
\widetilde{T}^{\varepsilon,\,3}_N&=&|(w_N^\varepsilon.\nabla
v_N^\varepsilon/w_N^\varepsilon)_{H^{0,s}}|,\\
\widetilde{T}^{\varepsilon,\,4}_N&=&|(v_N^\varepsilon.\nabla
v_N^\varepsilon/w_N^\varepsilon)_{H^{0,s}}|.\end{eqnarray*}Let us
now bound these four terms.\\As in
\cite{paicu'}, the following lemma is the main ingredient to estimate the first term.
We postpone the proof to the fifth section.
\begin{lem}
\label{u,v,0,s} Let $s>1/2$ be a real number. For any vector fields
$u$ and $v$ belonging to $\mathcal{H}_0^{0,s}(\Omega)\cap
\mathcal{H}_0^{1,s}(\Omega)$, we have
\begin{eqnarray*}
|(u.\nabla v/ v)_{H^{0,s}}|&\leq& C
(\|u\|_{H^{0,s}}^{\frac{1}{2}}\|\nabla_hu\|_{H^{0,s}}^{\frac{1}{2}}
\|v\|_{H^{0,s}}^{\frac{1}{2}}\|\nabla_hv\|_{H^{0,s}}^{\frac{3}{2}}\\&+&
\|\nabla_hu\|_{H^{0,s}}\|v\|_{H^{0,s}}\|\nabla_hv\|_{H^{0,s}}).\end{eqnarray*}
In particular, if $u=v$, then
$$
|(u.\nabla u/ u)_{H^{0,s}}|\leq C\|u\|_{H^{0,s}}\|\nabla_h
u\|_{H^{0,s}}^2.$$\end{lem}

\noindent Thanks to this lemma, we
have$$\widetilde{T}^{\varepsilon,\,1}_N\leq
C\|w_N^\varepsilon\|_{H^{0,s}}\|\nabla_h
w_N^\varepsilon\|_{H^{0,s}}^2.$$Hence, we get on the interval
$[0,T_{\varepsilon,N}[$
\begin{eqnarray*}
 \widetilde{T}^{\varepsilon,\,1}_N\leq& \frac{\nu_h}{100}\|\nabla_h
w_N^\varepsilon\|_{H^{0,s}}^2.\end{eqnarray*}Applying again Lemma
\ref{u,v,0,s} for the second term, we obtain
\begin{eqnarray*}\widetilde{T}^{\varepsilon,\,2}_N\leq &C&
\|v_N^\varepsilon\|_{H^{0,s}}^{\frac{1}{2}}\|\nabla_h
v_N^\varepsilon\|_{H^{0,s}}^{\frac{1}{2}}\|w_N^\varepsilon\|_{H^{0,s}}^{\frac{1}{2}}
\|\nabla_hw_N^\varepsilon\|_{H^{0,s}}^{\frac{3}{2}}\\&+&C\|\nabla_hv_N^\varepsilon\|_{H^{0,s}}\|
w_N^\varepsilon\|_{H^{0,s}}\|\nabla_hw_N^\varepsilon\|_{H^{0,s}}.\end{eqnarray*}
Using the convexity inequality $ab\leq \theta
a^{\frac{1}{\theta}}+(1-\theta)b^{\frac{1}{1-\theta}}$, we get on
the interval $[0,T_{\varepsilon,N}[$
\begin{eqnarray*}
\widetilde{T}^{\varepsilon,\,2}_N\leq
C\|w_N^\varepsilon\|_{H^{0,s}}^{2}
\|\nabla_hv_N^\varepsilon\|^{2}_{H^{0,s}}(1+\|
v_N^\varepsilon\|_{H^{0,s}}^2) +\frac{\nu_h}{50} \|\nabla_h
w_N^\varepsilon\|_{H^{0,s}}^{2}.
\end{eqnarray*}
Now we have to estimate the third term
$\widetilde{T}^{\varepsilon,\,3}_N$. First of all, we split it into
two parts $$\widetilde{T}^{\varepsilon,\,3}_N
=(w_{N,h}^\varepsilon.\nabla_h v_N^\varepsilon/
w_N^\varepsilon)_{H^{0,s}}+(w_{N,3}^\varepsilon.\partial_3
v_N^\varepsilon/ w_N^\varepsilon)_{H^{0,s}}.$$The Cauchy-Schwarz
inequality leads to $$|(w_{N,h}^\varepsilon.\nabla_h
v_N^\varepsilon/ w_N^\varepsilon)_{H^{0,s}}|\leq
\|w_{N,h}^\varepsilon.\nabla_h
v_N^\varepsilon\|_{H^{0,s}}\|w_N^\varepsilon\|_{H^{0,s}}.$$ The
following lemma will be useful for the estimate of
$\|w_{N,h}^\varepsilon.\nabla_h v_N^\varepsilon\|_{H^{0,s}}$. We
refer to the fifth section for the proof.
\begin{lem}\label{self}
\label{u,v,1,s}Let $s>\frac{1}{2}$ be a real number. For any
axisymmetric vector fields $u$ and $v$ such that $u\in\mathcal{H}_0^{0,s}(\Omega)\cap \mathcal{H}_0^{1,s}(\Omega)$ and
$v\in\mathcal{H}_0^{0,s}(\Omega)$, we have
$$\|u v\|_{H^{0,s}}\leq C\| v\|_{H^{0,s}}\left(\|u
\|_{H^{0,s}}+\|\nabla_hu \|_{H^{0,s}}\right).$$\end{lem} Hence, we get
$$\|w_{N,h}^\varepsilon.\nabla_h v_N^\varepsilon\|_{H^{0,s}}\leq C\|
\nabla_h
v_N^\varepsilon\|_{H^{0,s}}(\|w_N^\varepsilon\|_{H^{0,s}}+\|\nabla_hw_N^\varepsilon
\|_{H^{0,s}}),$$ and then, we
obtain\begin{eqnarray*}|(w_{N,h}^\varepsilon.\nabla_h
v_N^\varepsilon/ w_N^\varepsilon)_{H^{0,s}}|&\leq& C \|\nabla_h
v_N^\varepsilon\|_{H^{0,s}}
\|w_N^\varepsilon\|_{H^{0,s}}^2+C\|w_N^\varepsilon\|_{H^{0,s}}\|\nabla_hv_N^\varepsilon\|_{H^{0,s}}
\|\nabla_hw_N^\varepsilon\|_{H^{0,s}} \\
&\leq&\frac{\nu_h}{100}\|\nabla_h
w_N^\varepsilon\|_{H^{0,s}}^2+C\|w_N^\varepsilon\|_{H^{0,s}}^2(\|\nabla_hv_N^\varepsilon\|_{H^{0,s}}+
\|\nabla_hv_N^\varepsilon\|_{H^{0,s}}^2).\end{eqnarray*}Using again
the Cauchy-Schwarz inequality, we get
$$|(w_{N,3}^\varepsilon.\partial_3 v_N^\varepsilon/
w_N^\varepsilon)_{H^{0,s}}|\leq \|w_{N,3}^\varepsilon.\partial_3
v_N^\varepsilon\|_{H^{0,s}}\|w_N^\varepsilon\|_{H^{0,s}}.$$ Lemma
\ref{u,v,1,s} implies that
\begin{eqnarray*}\|w_{N,3}^\varepsilon.\partial_3
v_N^\varepsilon\|_{H^{0,s}}&\leq& C\|\partial_3
v_N^\varepsilon\|_{H^{0,s}}(\|
w_N^\varepsilon\|_{H^{0,s}}+\|\nabla_h
w_N^\varepsilon\|_{H^{0,s}}).\end{eqnarray*}To estimate the term
$\|\partial_3 v_N^\varepsilon\|_{H^{0,s}}$, it is worth noting that
$S_N^{x_3}v_N^\varepsilon=v_N^\varepsilon$. As $B$ is independent of
$x_3$, we have
$$S_N^{x_3}\Big(
v_N^\varepsilon\times B\Big)= (S_N^{x_3}v_N^\varepsilon)\times
B.$$Thus, $v_N^\varepsilon$ and $S_N^{x_3}v_N^\varepsilon$ satisfy
the same equation. Moreover, we notice that
$$v^\varepsilon_{N/t=0}=S_N^{x_3}v^\varepsilon_{N/t=0}=S_N^{x_3}u_0.$$
By uniqueness we get
$S_N^{x_3}v_N^\varepsilon=v_N^\varepsilon$. Consequently, we have
$$\| \partial_3 v_N^\varepsilon\|_{H^{0,s}}\leq C\| v_N^\varepsilon\|_{H^{0,s}},$$
and then$$|(w_{N,3}^\varepsilon.\partial_3 v_N^\varepsilon/
w_N^\varepsilon)_{H^{0,s}}|\leq \frac{\nu_h}{100} \|\nabla_h
w_N^\varepsilon\|_{H^{0,s}}^2+C\|w_N^\varepsilon\|_{H^{0,s}}^2(\|v_N^\varepsilon\|_{H^{0,s}}^2
+\|v_N^\varepsilon\|_{H^{0,s}}).$$ So, it turns out that
\begin{eqnarray*} \widetilde{T}^{\varepsilon,\,3}_N\leq
&C&\|w_N^\varepsilon\|_{H^{0,s}}^2(\|v_N^\varepsilon\|_{H^{0,s}}
+\|v_N^\varepsilon\|_{H^{0,s}}^2
+\|\nabla_hv_N^\varepsilon\|_{H^{0,s}}+\|\nabla_hv_N^\varepsilon\|_{H^{0,s}}^2)\\&+&\frac{\nu_h}{50}
\|\nabla_h w_N^\varepsilon\|_{H^{0,s}}^2.\end{eqnarray*}Now, we deal
with the fourth term $\widetilde{T}^{\varepsilon,\,4}_N$. As
$v_N^\varepsilon$ is a divergence free vector field , an integration
by parts leads to
\begin{eqnarray*}
\widetilde{T}^{\varepsilon,\,4}_N&=&(div(v_N^\varepsilon \otimes
v_N^\varepsilon)/ w_N^\varepsilon)_{H^{0,s}}\\&=&(v_N^\varepsilon
\otimes v_N^\varepsilon/\nabla_h
w_N^\varepsilon)_{H^{0,s}}+(\partial_3(v_N^\varepsilon \otimes
v_N^\varepsilon)/w_N^\varepsilon)_{H^{0,s}}.
\end{eqnarray*}

First, we estimate the term $(v_N^\varepsilon \otimes
v_N^\varepsilon/\nabla_h w_N^\varepsilon)_{H^{0,s}}$. Thanks to the
Cauchy-Schwarz inequality, we obtain
$$(v_N^\varepsilon \otimes v_N^\varepsilon/\nabla_h
w_N^\varepsilon)_{H^{0,s}}\leq \|v_N^\varepsilon \otimes
v_N^\varepsilon\|_{H^{0,s}} \|\nabla_h
w_N^\varepsilon\|_{H^{0,s}}.$$The fact that $S_N^{x_3}v
_N^\varepsilon =v_N^\varepsilon$ implies that
\begin{eqnarray*}\|v_N^\varepsilon \otimes
v_N^\varepsilon\|_{H^{0,s}}&\leq& C\|v_N^\varepsilon \otimes
v_N^\varepsilon\|_{L^2(\Omega)}\\&\leq&C\|v_N^\varepsilon \otimes
v_N^\varepsilon\|_{L^2_v( L^2_h)}.\end{eqnarray*}As
$\|v_N^\varepsilon \otimes v_N^\varepsilon\|_{L^2_v( L^2_h)}$ is
equivalent to $\| \widetilde{v_N^\varepsilon }\otimes
\widetilde{v_N^\varepsilon}\|_{L^2_v(L^2(rdr))}$, we get by the
H\"older inequality
$$\|\widetilde{v_N^\varepsilon} \otimes
\widetilde{v_N^\varepsilon}\|_{L^2_v( L^2(rdr))}\leq \|
\widetilde{v_N^\varepsilon}\|_{L^\infty_v(L^4(rdr))}\|
\widetilde{v_N^\varepsilon}\|_{L^2_v(L^4(rdr))}.$$Since we have
$S_N^{x_3}v_N^\varepsilon=v_N^\varepsilon$, Bernstein lemma yields
$$\|
\widetilde{v_N^\varepsilon}\|_{L^\infty_v(L^4(rdr))}\leq C \|
\widetilde{v_N^\varepsilon}\|_{L^2_v(L^4(rdr))}.$$ Hence, we get
$$\|v_N^\varepsilon \otimes v_N^\varepsilon\|_{L^2(\Omega)}\leq C \|
\widetilde{v_N^\varepsilon}\|_{L^2_v(L^4(rdr))}^2.$$ Notice that
$\|\widetilde{v_N^\varepsilon}(.,x_3)\|_{L^4(rdr)}=\|r^{\frac{1}{4}}
\widetilde{v_N^\varepsilon}(.,x_3)\|_{L^4(dr)}$. We suppose that
$x_3$ is fixed and we define
$\psi_{x_3}^{\varepsilon,N}(r):=r^{\frac{1}{4}}
\widetilde{v_N^\varepsilon}(r,x_3)$. Using the assumption made on
the domain $\Omega$, we infer
$$\|\psi_{x_3}^{\varepsilon,N}\|_{L^2(dr)}\leq
\frac{1}{\rho^{\frac{1}{4}}}\|\widetilde{v_N^\varepsilon}(.,x_3)\|_{L^2(rdr)}$$
and$$\|\partial_r\psi_{x_3}^{\varepsilon,N}\|_{L^2(dr)}\leq
\frac{C}{\rho^{\frac{5}{4}}}\|\widetilde{v_N^\varepsilon}(.,x_3)\|_{L^2(rdr)}+
\frac{C}{\rho^{\frac{1}{4}}}\|\partial_r\widetilde{v_N^\varepsilon}
(.,x_3)\|_{L^2(rdr)}.$$Since $v_N^{\varepsilon}$ is axisymmetric,
then $\|\nabla_h v_N^{\varepsilon}(.,x_3)\|_{L^2_h}$ and
$\|\partial_r \widetilde{v^\varepsilon_N}(.,x_3)\|_{L^2(rdr)}$ are
equivalent. But, as $v_N^{\varepsilon}(.,x_3)$ is in $L^2(\Omega_h)$
with $\nabla_h v_N^{\varepsilon}(.,x_3)$ in $L^2(\Omega_h)$, then we
get $\psi_{x_3}^{\varepsilon,N}\in H^1(]\rho,+\infty[)$.\\In order
to use Sobolev embedding on the whole space $\mathbb{R}$, we need to
extend the function $\psi_{x_3}^{\varepsilon,N}$. The
following extension lemma is needed.

\begin{lem}\label{equi}There exists an
extension operator
$$P:H^1(]\rho,+\infty[)\rightarrow H^1(\mathbb{R})$$such that
\begin{eqnarray*}Pu_{/]\rho,+\infty[}&=&u,\\
\|Pu\|_{L^2(\mathbb{R})}&\leq& C\|u\|_{L^2(]\rho,+\infty[)},\\
\|Pu\|_{H^1(\mathbb{R})}&\leq&
C\|u\|_{H^1(]\rho,+\infty[)}.\end{eqnarray*}\end{lem}

\begin{rem}
We can consider (for example)
\begin{eqnarray*}Pu(r,x_3)=
\begin{cases}
u(r,x_3),\,\,\, \mbox{if}\,\,\,
r>\rho\\u(2\rho-r,x_3),\,\,\,\mbox{if}\,\,\,{r\leq\rho}.
\end{cases}
\end{eqnarray*}More general results about extension operators can be found in \cite{Bresis} for instance.
\end{rem}

\smallskip

Now, we continue the study of the term $(v_N^\varepsilon \otimes
v_N^\varepsilon/\nabla_h w_N^\varepsilon)_{H^{0,s}}$. Thanks to Sobolev embeddings, we get
$$\|P\psi_{x_3}^{\varepsilon,N}\|_{L^4(dr)}\leq C
\|P\psi_{x_3}^{\varepsilon,N}\|_{\dot{H}^\frac{1}{4}(\mathbb{R})}.$$By
interpolation, we
have\begin{eqnarray*}\|P\psi_{x_3}^{\varepsilon,N}\|_{\dot{H}^\frac{1}{4}(\mathbb{R})}&\leq&
\|\partial_r (P\psi_{x_3}^{\varepsilon,N})\|_{L^2(dr)}^{\frac{1}{4}}
\|P\psi_{x_3}^{\varepsilon,N}\|_{L^2(dr)}^{\frac{3}{4}}\\&\leq&
\|\psi_{x_3}^{\varepsilon,N}\|_{H^1(]\rho,+\infty[)}^{\frac{1}{4}}\|\psi_{x_3}^{\varepsilon,N}\|_
{L^2(]\rho,+\infty[)}^{\frac{3}{4}}
\\&\leq& \|\psi_{x_3}^{\varepsilon,N}\|_{L^2(dr)}+\|\psi_{x_3}^{\varepsilon,N}\|_{L^2(dr)}^{\frac{3}{4}}
\|\partial_r\psi_{x_3}^{\varepsilon,N}\|_{L^2(dr)}^{\frac{1}{4}}.\end{eqnarray*}
Hence, we get
$$\|\widetilde{v_N^\varepsilon}(.,x_3)\|_{L^4(rdr)}\leq
C\|\widetilde{v_N^\varepsilon}(.,x_3)\|_{L^2(rdr)}+C\|\widetilde{v_N^\varepsilon}(.,x_3)\|_{L^2(rdr)}
^{\frac{3}{4}}\|\partial_r\widetilde{v_N^\varepsilon}(.,x_3)\|_{L^2(rdr)}
^{\frac{1}{4}}.$$ Taking the $L^2_v$ norm, we obtain thanks to the
H\"older inequality$$\|
\widetilde{v_N^\varepsilon}\|_{L^2_v(L^4(rdr))}\leq C \|\partial_r
\widetilde{v_N^\varepsilon}\|_{L^2_v (L^2(rdr))}^{\frac{1}{4}}
\|\widetilde{v_N^\varepsilon}\|_{L^2 _v
(L^2(rdr))}^{\frac{3}{4}}+C\|\widetilde{v_N^\varepsilon}\|_{L^2 _v
(L^2(rdr))},$$and then
$$(v_N^\varepsilon \otimes
v_N^\varepsilon/\nabla_h w_N^\varepsilon)_{H^{0,s}}\leq
\frac{\nu_h}{100}\|\nabla_hw_N^\varepsilon\|_{H^{0,s}}^2+C\|\nabla_hv_N^\varepsilon\|_{L^2(\Omega)}\|v_N^\varepsilon
\|_{L^2(\Omega)}^{3}+C \|v_N^\varepsilon\|_{L^2(\Omega)}^{4}.$$ As
$v_N^\varepsilon$ is localized in vertical frequencies, the term
$(\partial_3(v_N^\varepsilon \otimes
v_N^\varepsilon)/w_N^\varepsilon)_{H^{0,s}}$ is estimated as the
term $(v_N^\varepsilon \otimes v_N^\varepsilon/\nabla_h
w_N^\varepsilon)_{H^{0,s}}$ and we get $$(\partial_3(v_N^\varepsilon
\otimes v_N^\varepsilon)/w_N^\varepsilon)_{H^{0,s}}\leq C
\|w_N^\varepsilon\|_{H^{0,s}}^2+C\|\nabla_hv_N^\varepsilon\|_{L^2(\Omega)}\|v\|_{L^2(\Omega)}^{3}+
C\|v_N^\varepsilon\|_{L^2(\Omega)}^{4}.$$Considering all the above
estimates, we get on the interval $[0,T_{\varepsilon ,N}[$
\begin{eqnarray*}
\frac{1}{2}\frac{d}{dt}\|w_N^\varepsilon\|_{H^{0,s}}^2+\nu_h
\|\nabla_h w_N^\varepsilon\|_{H^{0,s}}^2& \leq& C
\|w_N^\varepsilon\|_{H^{0,s}}^2(1+\|v_N^\varepsilon\|_{H^{0,s}}+
\|v_N^\varepsilon\|_{H^{0,s}}^2+\|\nabla_h
v_N^\varepsilon\|_{H^{0,s}}\\&+&\|\nabla_hv_N^\varepsilon\|_{H^{0,s}}^2+\|v_N^\varepsilon\|_{H^{0,s}}^{2}
\|\nabla_hv_N^\varepsilon\|_{H^{0,s}}^2 )+ \frac{3}{50}\nu_h
\|\nabla _hw_N^\varepsilon\|_{H^{0,s}}^2         \\&+&C(
\|v_N^\varepsilon\|_{L^2(\Omega)}^3\|\nabla_h
v_N^\varepsilon\|_{L^2(\Omega)}+
\|v_N^\varepsilon\|_{L^2(\Omega)}^4).
\end{eqnarray*}
By using Gronwall's lemma, we deduce that
\begin{eqnarray*}
\|w^{\varepsilon}_N (t)\|^2_{H^{0,s}}&\leq& \Big[\|w^{\varepsilon}_N
(0)\|_{H^{0,s}}^2+\int_0^t
C(\|v_N^\varepsilon(\tau)\|_{L^2(\Omega)}^3\|\nabla_h
v_N^\varepsilon(\tau)\|_{L^2(\Omega)}+\|v_N^\varepsilon(\tau)\|_{L^2(\Omega)}^4)d\tau\Big]\\
&\times&exp\bigg[\int_0^t C\big(1+\|
v_N^\varepsilon(\tau)\|_{H^{0,s}}+\|v_N^\varepsilon(\tau)\|_{H^{0,s}}^2+
\|\nabla_hv_N^\varepsilon(\tau)\|_{H^{0,s}}+\|\nabla_h
v_N^\varepsilon(\tau)\|_{H^{0,s}}^2\\&&+\|v_N^\varepsilon(\tau)\|_{H^{0,s}}^2\|\nabla_h
v_N^\varepsilon(\tau)\|_{H^{0,s}}^{2}
 \big)d\tau\bigg]
.\end{eqnarray*} It is of interest to note that $L^2$ and $H^{0,s}$
energy estimates on $v^\varepsilon_N$ imply that
$$
\|v^\varepsilon_N(t)\|^2_{L^2(\Omega)}+2\nu_h\int_0^t\|\nabla_h
v^\varepsilon_N(\tau)\|^2_{L^2(\Omega)}d\tau \leq
\|u_0\|_{L^2(\Omega)}^2,
$$
and
$$
\|v^\varepsilon_N(t)\|^2_{H^{0,s}}+2\nu_h\int_0^t\|\nabla_h
v^\varepsilon_N(\tau)\|^2_{H^{0,s}}d\tau\leq \|u_0\|_{H^{0,s}}^2.
$$Thanks to the H\"older inequality, we obtain
$$\int_0^t C(\|v_N^\varepsilon(\tau)\|_{L^2(\Omega)}^3\|\nabla_h
v_N^\varepsilon(\tau)\|_{L^2(\Omega)}+\|v_N^\varepsilon(\tau)\|_{L^2(\Omega)}^4)d\tau\leq
C\|u_0\|_{L^2(\Omega)}^4( t+\sqrt{t}).$$Finally, we
get$$\|w^{\varepsilon}_N (t)\|^2_{H^{0,s}}\leq
\big[\|w^{\varepsilon}_N
(0)\|_{H^{0,s}}^2+C\|u_0\|_{L^2(\Omega)}^4(t+\sqrt{t})\big]\times
exp(C+Ct+C\sqrt{t}),
$$where $C>0$ depends on $\|u_0\|_{H^{0,s}}$.\\Let us consider a positive real number $T$ such that
$$\big[\|w^{\varepsilon}_N
(0)\|_{H^{0,s}}^2+C\|u_0\|_{L^2(\Omega)}^4(T+\sqrt{T})\big]\times
exp(C+CT+C\sqrt{T})\leq (\frac{3}{2}c\nu_h)^2.
$$Notice that $T$ is independent of $\varepsilon$. Therefore,
$w^{\varepsilon}_N$ exists on the time interval [0,T]. But as
$v^\varepsilon_N$ is global in time, then
$u^\varepsilon=v^\varepsilon_N+w^\varepsilon_N$ exists on the time
interval [0,T] and Theorem \ref{th5} is proved.\end{proof}As said
before, uniform local existence in isotropic Sobolev space is a
consequence of Theorem \ref{th5}. We note that the proof of this
fact is contained in \cite{paicu'}, but we give it for the
convenience of the reader.
\begin{proof}[Proof of Theorem \ref{th6}]\quad\\
As $u_0$ belongs to the space $\mathcal{H}^m_0(\Omega)$, $m\geq 1$,
then $u_0$ is in $\mathcal{H}^{0,m}_0(\Omega)$. Hence, Theorem
\ref{th5} yields the existence of a unique solution $u^\varepsilon$
for the system $(S^\varepsilon)$ on a uniform time interval $[0,T]$
such that $\|u^\varepsilon(t)\|_{H^{0,s}}$ and $\dint_0^t\|\nabla_h
u^\varepsilon(\tau)\|_{H^{0,s}}^2d\tau$ are uniformly bounded on the
time interval $[0,T]$. Let $T^\varepsilon$ be the maximal time of
existence of $u^\varepsilon$ in $\mathcal{H}^m_0(\Omega)$. Theorem
\ref{th1} and the inclusion of $\mathcal{H}_0^{0,m}(\Omega)$ in
$L^\infty_v(L^2_h)$ imply that $T^\varepsilon>T$. Thus, we get the
uniform local in time existence of $u^{\varepsilon}$ in the space
$\mathcal{H}^m_0(\Omega)$.\end{proof}
%@@@@@@@@@@@@@@@@@@@@@@@@@@@@@@@@@@@@%@@@@@@@@@@@@@@@@@@@@@@@@@@@@@@@@@@@@%@@@@@@@@@@@@@@@
\section{\label{final}Product laws}
Before proving the technical lemmas, let us first recall some
results about the anisotropic Littlewood Paley
theory.\subsection{Anisotropic Littlewood Paley theory} Anisotropic
Sobolev spaces can be characterized using a dyadic decomposition in
the vertical frequency space. So, let us first recall some elements
of the Littlewood-Paley theory, the details of which can be found in
\cite{Iftimie} for instance.\\Let $u$ be a function defined on $\Omega$, we
have$$\mathcal{F}^{\vee}(\Delta_{q}^{\vee} u)(x_h,.) =
 \varphi(\frac{|.|}{2^{q}})\mathcal{F}^{\vee}(u)(x_h,.), \quad q \geq 0,$$
 $$\mathcal{F}^{\vee}(\Delta_{-1}^{\vee}u)(x_h,.)=
\chi(|.|)\mathcal{F}^{\vee}(u)(x_h,.) ,$$
 $$\Delta_{q}^{\vee} u = 0 ,\quad q \leq - 2.$$
 The positive functions $\varphi$ and $\chi$ represent a dyadic
 partition of unity in $\mathbb{R}$, that is to say they are smooth
 functions such that
 $$supp \,\chi \subset B(0, \frac{4}{3}), \,supp\, \varphi \subset \mathcal{C}
 (0, \frac{3}{4}, \frac{8}{3}), $$ and $\forall t \in \mathbb{R}$
 $$\chi(t) + \sum_{q \geq 0} \varphi(2^{-q} t) = 1.$$
Let us also define the operator
$$S^\vee_q u = \sum_{q' \leq q-1}\Delta^\vee_{q'}u.$$ Those
definitions enable us to characterize  anisotropic Sobolev spaces
$H^{0,s}(\Omega)$. More precisely, a tempered distribution $u$
belongs to $H^{0,s}(\Omega)$ if and only if $$\sum_q 2^{2qs}
\|\Delta^{\vee}_q u\|^2_{L^2(\Omega)}<\infty.$$ Moreover, we have
$$ \|u\|^2_ {H^{0,s}(\Omega)}\approx\sum_q 2^{2qs} \|\Delta^{\vee}_q u\|^2_{L^2(\Omega)}.$$
 The dyadic decomposition is
also important for studying the product of two distributions thanks
to Bony's decomposition.\\
Let $u$ and $v$ be two distributions. We have$$u = \sum_{q \in
\mathbb{Z}}\Delta^\vee_q u \, ; \, v = \sum_{q \in \mathbb{Z}}
\Delta^\vee_q v.$$ We denote \begin{eqnarray*}T_u
v&=&\sum_{q}S^\vee_{q-1} u.\Delta^\vee_q v  \\
R(u,v)&=&\sum_{{\tiny{\begin{array}{c}q\\i\in\{0,\pm 1\}\\
\end{array}}}}     \Delta^\vee_q u.\Delta^\vee_{q-i} v \,.\end{eqnarray*}
Thus, we obtain
$$u.v = T_u v+T_v u+R(u,v).$$
We have
\begin{eqnarray*}
\Delta^\vee_q(uv)&=&\sum_{|q'-q|\leq
4}\Delta_q^\vee(S^\vee_{q'-1}u.\Delta_{q'}^\vee v)+ \sum_{|q'-q|\leq
4}\Delta_q^\vee(S^\vee_{q'-1}v
.\Delta_{q'}^\vee u)\\&+&\sum_{{\tiny{\begin{array}{c}i\in\{0,\pm 1\}\\  q'>q-4\\
\end{array}}}}\Delta_q^\vee(\Delta_{q'}^\vee u.\Delta_{q'-i}^\vee
v).
\end{eqnarray*}
The dyadic decomposition is useful in the sense that the derivatives
in vertical variable act in a very special way on functions
localized in vertical frequencies in a ball or a ring. More
precisely, we have the following lemma the proof of which can be
found in \cite{paicu'}, \cite{paicu} for instance.
\begin{lem}{\bf{Bernstein lemma}}\quad\\
Let $p, r$ and $r'$be numbers such that $\infty \geq p \geq 1$ and
$\infty\geq r \geq r' \geq 1$.\\ Then, there exists a constant $C>0$
such that for any vector field $u $ defined on
$\Omega_h\times\mathbb{R}$ with supp$\mathcal{F}^\vee u \subset
\mathbb{R}^2_h \times 2^q \mathcal{C}$, where $\mathcal{C}$ is a
dyadic ring, we have
$$2^{qk} C^{-k}\|u\|_{L^p_h(L^r_v)}\, \leq\,
\|\partial^k_{x_3}u\|_{L^p_h(L^r_{v})} \,\leq\, 2^{qk} C^k
\|u\|_{L^p_h(L^r_{v})},$$
$$2^{qk}C^{-k}\|u\|_{L^{r}_v(L^p_h)} \leq \|\partial^k_{x_3}u\|_{L^r_v(L^p_h)}
\leq 2^{qk}C^k \|u\|_{L^r_v(L^p_h)},$$
$$\|u\|_{L^p_h(L^r_v)} \leq C 2^{q(\frac{1}{r'} - \frac{1}{r})}
\|u\|_{L^p_h(L^{r'}_v)},$$
$$\|u\|_{L^r_v(L^p_h)} \leq C 2^{q(\frac{1}{r'} - \frac{1}{r})}
\|u\|_{L^{r'}_v(L^p_h)}.$$\end{lem}
\subsection{Proofs of the technical lemmas}
In this part, we denote by $(b_q)_{q\in \mathbb{Z}}$ and
$(c_q)_{q\in \mathbb{Z}}$ positive sequences such that
$$\sum_{q\in \mathbb{Z}} b_{q} \leq 1\quad \mbox{and}\quad \sum_{q\in \mathbb{Z}}c_q^{2}\,\leq\,1.$$
For the proof of Lemma \ref{u,v,0,s}, we proceed as in \cite{paicu}
where the critical Besov space ${\mathcal B}^{0,\frac{1}{2}}$ is used.
\begin{proof}[Proof of Lemma \ref{u,v,0,s}]\quad\\
The proof of Lemma \ref{u,v,0,s} relies on basic inequalities.
\begin{prop}\label{DP}\quad\\
For any vector field $u$ in $H^{0,s}(\Omega)$, we have
$$(1)\quad \|\Delta^\vee_{q}u\|_{L^2}\,\leq\,C \,c_q\,
2^{-qs}\|u\|_{H^{0,s}},$$
$$(2)\quad \|u\|_{L^{\infty}_{v}(L^2_h)}\leq C \|u\|_{H^{0,s}}, \,\,s>\frac{1}{2}.$$
For any vector field $u$ in $\mathcal{H}_0^{0,s}(\Omega)\cap
\mathcal{H}_0^{1,s}(\Omega)$, we have
$$(3)\quad \|\Delta^\vee_{q}u\|_{L^2_v
(L^4_h)}\leq C  \,c_q\,
2^{-qs}\|u\|_{H^{0,s}}^{\frac{1}{2}}\|\nabla_h
u\|_{H^{0,s}}^{\frac{1}{2}},$$
$$(4)\quad\|u\|_{L^{\infty}_v (L^4_h)}\leq C
\|u\|_{H^{0,s}}^{\frac{1}{2}}\|\nabla_h
u\|_{H^{0,s}}^{\frac{1}{2}},\,\, s>\frac{1}{2}.$$
\end{prop}
\begin{proof}[Proof of Proposition \ref{DP}]\quad\\
Thanks to the Bernstein lemma, we get
$$\|u\|_{L^{\infty}_{v}(L^2_h)}\leq C \sum_q 2^{\frac{q}{2}}\|\Delta^\vee_{q}u\|_{L^2}.$$
As $s>\frac{1}{2}$, then the Cauchy-Schwarz inequality leads to (2).\\
To get (3), we just have to prove it for $u$ in
$C^{\infty}_c(\Omega)$. Sobolev embeddings imply that
\begin{eqnarray*}\|\Delta_q^\vee u(.,x_3)\|_{L^4_h}&\leq& C
\|\Delta_q^\vee
u(.,x_3)\|_{\dot{H}^{\frac{1}{2}}(\mathbb{R}^2)}\\&\leq&
C\|\Delta_q^\vee u(.,x_3)\|_{L^2_h}^{\frac{1}{2}}
\|\nabla_h\Delta_q^\vee
u(.,x_3)\|_{L^2_h}^{\frac{1}{2}}.\end{eqnarray*}Taking the $L^2_v$
norm and using (1), we get (3).\\To prove $(4)$, we use the
Bernstein lemma to get
$$\|u\|_{L^{\infty}_v (L^4_h)}\leq C \sum_q 2^{\frac{q}{2}}
\|\Delta^\vee_{q}u\|_{L^2_{v}(L^4_h)}.$$ As $s>\frac{1}{2}$, we get
the result thanks to (3) and to the Cauchy-Schwarz inequality.
\end{proof}
Let us go back to the proof of Lemma \ref{u,v,0,s}.\\

We have \begin{eqnarray*} (u.\nabla
v/v)_{H^{0,s}}&=&\sum_q2^{2qs}(\Delta_q^\vee (u.\nabla v)/
\Delta_q^\vee v)_{L^2}\\&=&\sum_q 2^{2qs}(F_q^h/\Delta_q^\vee
v)_{L^2}+\sum_q 2^{2qs}(F_q^\vee /\Delta_q^\vee v)_{L^2},
\end{eqnarray*}
where $F_q^h=\Delta_q^\vee (u_h.\nabla_h v)$ and
$F_q^\vee=\Delta_q^\vee (u_3.\partial_3 v)$.\\For the term
$(F_q^h/\Delta_q^\vee v)_{L^2}$, we have thanks to the H\"older
inequality$$|(F_q^h/\Delta_q^\vee v)_{L^2}|\leq \|F_q^h\|_{L^2_v
(L^{\frac{4}{3}}_h) }\|\Delta_q^\vee v\|_{L^2_v
(L^4_h)}.$$Proposition \ref{DP} leads to
\begin{eqnarray*} \|\Delta_q^\vee v\|_{L^2_v (L^4_h)}\leq
C  \,c_q\, 2^{-qs}\|v\|_{H^{0,s}}^{\frac{1}{2}}\|\nabla_h
v\|_{H^{0,s}}^{\frac{1}{2}}.
\end{eqnarray*}Bony's decomposition implies
that$$F_q^h=T_{1,q}^h+T_{2,q}^h+T_{3,q}^h,$$where
\begin{eqnarray*}
T_{1,q}^h&=&\sum_{|q'-q|\leq
4}\Delta_q^\vee(S^\vee_{q'-1}u_h.\Delta_{q'}^\vee \nabla_h
v),\\T_{2,q}^h&=& \sum_{|q'-q|\leq
4}\Delta_q^\vee(S^\vee_{q'-1}\nabla_h v
.\Delta_{q'}^\vee u_h),\\T_{3,q}^h&=&\sum_{{\tiny{\begin{array}{c}i\in\{0,\pm 1\}\\  q'>q-4\\
\end{array}}}}\Delta_q^\vee(\Delta_{q'}^\vee u_h.\Delta_{q'-i}^\vee \nabla_h
v).
\end{eqnarray*}
The H\"older inequality leads to
$$\|T_{1,q}^h\|_{L^2_v (L_h^{\frac{4}{3}})}\leq \sum_{|q'-q|\leq
4}\|S^\vee_{q'-1}u\|_{L^\infty_v (L^4_h)}\|\Delta_{q'}^\vee \nabla_h
v\|_{L^2}.$$But, we have
$$\|S^\vee_{q'-1}u\|_{L^\infty_v (L^4_h)}\leq C\|u\|_{L^{\infty}_v
(L^4_h)},$$then, by Proposition \ref{DP}, we obtain
\begin{eqnarray*}\|T_{1,q}^h\|_{L^2_v (L_h^{\frac{4}{3}})}&\leq&C 2^{-qs} \sum_{|q-q'|\leq4}c_{q'}2^{(q-q')s}
\|u\|_{H^{0,s}}^{\frac{1}{2}}\|\nabla_h
u\|_{H^{0,s}}^{\frac{1}{2}}\|\nabla_h v\|_{H^{0,s}} \\&\leq& C
\,c_q\, 2^{-qs}\|\nabla_h
v\|_{H^{0,s}}\|u\|_{H^{0,s}}^{\frac{1}{2}}\|\nabla_h
u\|_{H^{0,s}}^{\frac{1}{2}}.\end{eqnarray*} As for the term
$T_{2,q}^h$, the H\"older inequality implies that
\begin{eqnarray*}\|T_{2,q}^h\|_{L^2_v (L_h^{\frac{4}{3}})}&\leq&
\sum_{|q'-q|\leq 4}\|S^\vee_{q'-1}\nabla_h v\|_{L^\infty_v
(L^2_h)}\|\Delta_{q'}^\vee u\|_{L^2_{v}(L^4_h)}\\&\leq&
\|\nabla_hv\|_{L^\infty_v (L^2_h)}\sum_{|q'-q|\leq
4}\|\Delta_{q'}^\vee u\|_{L^2_{v}(L^4_h)}.\end{eqnarray*} Thanks to
Proposition \ref{DP}, we get$$\|T_{2,q}^h\|_{L^2_v
(L_h^{\frac{4}{3}})}\leq C \,c_q\, 2^{-qs}\|\nabla_h
v\|_{H^{0,s}}\|u\|_{H^{0,s}}^{\frac{1}{2}}\|\nabla_h
u\|_{H^{0,s}}^{\frac{1}{2}}.$$ For the term $T_{3,q}^h$, the
Bernstein lemma yields $$\|T_{3,q}^h\|_{L^2_v
(L_h^{\frac{4}{3}})}\leq
C 2^{\frac{q}{2}}\sum_{{\tiny{\begin{array}{c}i\in\{0,\pm 1\}\\  q'>q-4\\
\end{array}}}}\|\Delta_{q'}^\vee u_h.\Delta_{q'-i}^\vee \nabla_h
v\|_{L^1_v (L_h^{\frac{4}{3}})}.$$ By the H\"older inequality, we
obtain$$\|T_{3,q}^h\|_{L^2_v (L_h^{\frac{4}{3}})}\leq
C 2^{\frac{q}{2}}\sum_{{\tiny{\begin{array}{c}i\in\{0,\pm 1\}\\  q'>q-4\\
\end{array}}}}\|\Delta_{q'}^\vee u\|_{L^2_v (L^4_h)}\|\Delta_{q'-i}^\vee \nabla_h
v\|_{L^2}.$$Thanks to Proposition \ref{DP}, we
have$$\|T_{3,q}^h\|_{L^2_v (L_h^{\frac{4}{3}})}\leq C  \,c_q\,
2^{-qs}\|\nabla_h
v\|_{H^{0,s}}\|u\|_{H^{0,s}}^{\frac{1}{2}}\|\nabla_h
u\|_{H^{0,s}}^{\frac{1}{2}}.$$For the term $(F_q^\vee /\Delta_q^\vee
v)_{L^2}$, we use the following decomposition
$$(\Delta_q^\vee(u_3.\partial_3
v)/\Delta_q^\vee
v)_{L^2}=T_{1,q}^\vee+T_{2,q}^\vee+T_{3,q}^\vee+T_{4,q}^\vee,$$
where
\begin{eqnarray*}
T_{1,q}^\vee&=&\int_{\Omega_h\times\mathbb{R}}S^\vee_{q-1}u_3.\partial_3\Delta_q^\vee
v.\Delta_q^\vee v
\,dx,\\T_{2,q}^\vee&=&\sum_{|q-q'|\leq4}\int_{\Omega_h\times\mathbb{R}}[\Delta_q^\vee;S^\vee_{q'-1}u_3]
\partial_3\Delta_{q'}^\vee
v.\Delta_q^\vee v
\,dx,\\T_{3,q}^\vee&=&\sum_{|q-q'|\leq4}\int_{\Omega_h\times\mathbb{R}}(S^\vee_{q-1}u_3-S^\vee_{q'-1}u_3)
.\partial_3\Delta_q^\vee\Delta_{q'}^\vee v.\Delta_q^\vee v
\,dx,\\T_{4,q}^\vee&=&\sum_{q'>q-4}\int_{\Omega_h\times\mathbb{R}}\Delta_q^\vee(S^\vee_{q'+1}(\partial_3
v).\Delta_{q'}^\vee u_3).\Delta_q^\vee v \,dx.
\end{eqnarray*}
As $u$ is divergence free, we get after integration by
parts$$T_{1,q}^\vee=\frac{1}{2}\int_{\Omega_h\times\mathbb{R}}S^\vee_{q-1}(\textmd{div}_h
u_h)\Delta_q^\vee v.\Delta_q^\vee v\,dx.$$ The H\"older inequality
leads to\begin{eqnarray*}|T_{1,q}^\vee|&\leq&
C\|S^\vee_{q-1}(\textmd{div} _h u_h)\|_{L^\infty_v (L^2_h)}
\|\Delta_q^\vee v\|_{L^2_v (L^4_h)} ^2\\&\leq&C \|\nabla_h
u\|_{L^\infty_v (L^2_h)} \|\Delta_q^\vee v\|_{L^2_v (L^4_h)}^2.
\end{eqnarray*}Thanks to Proposition \ref{DP}, we
obtain$$|T_{1,q}^\vee| \leq C 2^{-2 q s}b_q \|\nabla_h
u\|_{H^{0,s}}\|v\|_{H^{0,s}}\|\nabla_h v\|_{H^{0,s}}.$$ For the term
$T_{2,q}^\vee$, we use the following lemma (see \cite {paicu} for
the proof).
\begin{lem}
Let $p, r, s$ and $t$ be real numbers such that $$1\leq p, r,s, t\leq\infty\quad\mbox{and}\quad
\displaystyle{\frac{1}{p} = \frac{1}{r} + \frac{1}{s}}.$$For any
vector fields $u$ and $v$ the following inequality holds
$$
\|[\Delta_q^{\vee}; u]v\|_{L^t_v L^p_h}\leq C 2^{-q} \|\partial_3
u\|_{L^\infty_v L^r_h} \|v\|_{L^t_v L^s _h},
$$
where $[\Delta_q^{\vee}; u]v=\Delta_q^{\vee}(u v)-u\Delta_q
^{\vee}v$.
\end{lem}
Hence, we get thanks to the H\"older inequality
\begin{eqnarray*}
|T_{2,q}^\vee|&\leq&
\sum_{|q-q'|\leq4}\|[\Delta_q^\vee;S^\vee_{q'-1}u_3]\partial_3\Delta_{q'}^\vee
v\|_{L^2_v (L^{\frac{4}{3}}_h)}\|\Delta_q^\vee v\|_{L^2_v
(L^4_h)}\\&\leq&C \sum_{|q-q'|\leq4}2^{-q}\|S^\vee_{q'-1}\partial_3
u_3\|_{L^\infty_v (L^{2}_h)}\|\partial_3\Delta_{q'}^\vee v\|_{L^2_v
(L^4_h)}\|\Delta_{q}^\vee v\|_{L^2_v (L^4_h)}.
\end{eqnarray*}
Since $ \Delta_{q'}^\vee v$ is localized in vertical frequencies in
the ring of size $2^{q'}$, we get
$$\|\partial_3\Delta_{q'}^\vee
v\|_{L^2_v (L^4_h)}\leq C 2^{q'}\|\Delta_{q'}^\vee v\|_{L^2_v
(L^4_h)}.$$As $\partial_3 u_3=-\textmd{div} _h u_h$, we obtain
thanks to Proposition \ref{DP}
$$|T^\vee_{2,q}|\leq C 2^{-2 q s}b_q \|\nabla_h
u\|_{H^{0,s}}\|v\|_{H^{0,s}}\|\nabla_h v\|_{H^{0,s}}.$$ As for the
term $T_{3,q}^\vee$, the H\"older inequality leads to
$$|T_{3,q}^\vee|\leq
\sum_{|q-q'|\leq4}\|S^\vee_{q'-1}u_3-S^\vee_{q-1}u_3\|_{L^{\infty}_{v}(L^2_h)}\|\partial_3\Delta_q^\vee
\Delta_{q'}^\vee v\|_{L^2_v(L^4_h)}.$$But, since
$(S^\vee_{q'-1}u_3-S^\vee_{q-1}u_3)$ is localized in vertical
frequencies in the ring of size $2^q$, we get thanks to the
Bernstein lemma \begin{eqnarray*}
\|S^\vee_{q'-1}u_3-S^\vee_{q-1}u_3\|_{L^{\infty}_{v}(L^2_h)}\leq C
2^{-q}\|(S^\vee_{q'-1}-S^\vee_{q-1})\partial_3u_3\|_{L^{\infty}_{v}(L^2_h)}.
\end{eqnarray*}
As $u$ is divergence free, we obtain by Proposition \ref{DP}
\begin{eqnarray*}\|(S^\vee_{q'-1}-S^\vee_{q-1})\partial_3u_3\|_{L^{\infty}_{v}(L^2_h)}\leq
C2^{-q} \|\nabla_h u\|_{H^{0,s}}.\end{eqnarray*}Applying again
Bernstein lemma and Proposition \ref{DP}, we get as in the above
estimates
$$|T_{3,q}^\vee|\leq C 2^{-2 q s}b_q \|\nabla_h
u\|_{H^{0,s}}\|v\|_{H^{0,s}}\|\nabla_h v\|_{H^{0,s}}.$$ Finally for
the term $T_{4,q}^\vee$, the H\"older inequality implies
that$$|T_{4,q}^\vee| \leq\sum_{q'>q-4}\|S^\vee_{q'+1}(\partial_3
v)\|_{L^\infty_v (L^4_h)}\|\Delta_{q'}^\vee
u_3\|_{L^2}\|\Delta_q^\vee v\|_{L^2_v (L^4_h)}.$$ Thanks to the
Bernstein lemma, we get$$\|S^\vee_{q'+1}(\partial_3 v)\|_{L^\infty_v
(L^4_h)}\leq C 2^{q'}\|v\|_{L^\infty_v (L^4_h)}.$$ Proposition
\ref{DP} leads to
$$\|S^\vee_{q'+1}(\partial_3 v)\|_{L^\infty_v (L^4_h)}\leq C
2^{q'}\|v\|_{H^{0,s}}^{\frac{1}{2}}\|\nabla_h
v\|_{H^{0,s}}^{\frac{1}{2}}.$$ Using again the Bernstein lemma and
the fact that $u$ is divergence free, we infer$$\|\Delta_{q'}^\vee
u_3\|_{L^2}\leq C 2^{-q'} \|\Delta_{q'}^\vee \nabla_h u\|_{L^2}.$$
We get thanks to Proposition \ref{DP}
$$|T_{4,q}^\vee| \leq C 2^{-2 q s}b_q \|\nabla_h
u\|_{H^{0,s}}\|v\|_{H^{0,s}}\|\nabla_h v\|_{H^{0,s}}.$$ Finally, we
sum all the estimates to conclude the proof of Lemma \ref{u,v,0,s}.
\end{proof}
\begin{proof}[Proof of Lemma \ref{u,v,1,s}]\quad\\
The proof of Lemma \ref{u,v,1,s} relies on some basic inequalities.
\begin{prop}\label{vic}
For any axisymmetric vector fields $u$ and $v$ such that $u$ belongs
to $\mathcal{H}^{0,s}_0(\Omega)\cap \mathcal{H}^{1,s}_0(\Omega)$ and
$v$ belongs to $\mathcal{H}^{0,s}_0(\Omega)$, we have
$$\|u(.,x_3)v(.,x_3)\|_{L^2_h}\leq C
\|v(.,x_3)\|_{L^2_h}(\|u(.,x_3)\|_{L^2_h}+\|\nabla_hu(.,x_3)\|_{L^2_h}).$$\end{prop}
\begin{proof}[Proof of Proposition \ref{vic}]\quad\\
By density arguments, we may suppose that $u$ are $v$ are in
$C^{\infty}_c(\Omega)$.\\We have
\begin{eqnarray*}
\|u(.,x_3)v(.,x_3)\|_{L^2_h}^2\leq
C\int_{\rho}^{+\infty}\tilde{u}^2(r,x_3)\tilde{v}^2(r,x_3)r
\,dr.\end{eqnarray*} But, we have
\begin{eqnarray*}
\tilde{u}^2(r,x_3)&=&\int_{\rho}^{r}2\tilde{u}(\sigma,x_3)\partial_r\tilde{u}(\sigma,x_3)\,d
\sigma\\&\leq&2\Big(\int_{\rho}^{r}|\tilde{u}(\sigma,x_3)|^2\,d
\sigma\Big)^{\frac{1}{2}}\Big(\int_{\rho}^{r}|\partial_r\tilde{u}(\sigma,x_3)|^2d
\sigma\Big)^{\frac{1}{2}}\\&\leq&2
\Big(\int_{\rho}^{+\infty}|\tilde{u}(\sigma,x_3)|^2\frac{1}{\sigma}\sigma
\,d\sigma\Big)^{\frac{1}{2}}\Big(\int_{\rho}^{+\infty}|\partial_r\tilde{u}(\sigma,x_3)|^2\frac{1}{\sigma}\sigma
\,d \sigma\Big)^{\frac{1}{2}}.\end{eqnarray*} Thus, we get
$$\tilde{u}^2(r,x_3)\leq C
\|\tilde{u}(.,x_3)\|_{L^2(rdr)}\|\partial_r\tilde{u}(.,x_3)\|_{L^2(rdr)},$$
and then$$\|u(.x_3)v(.,x_3)\|_{L^2_h}\leq C
\|v(.x_3)\|_{L^2_h}(\|u(.x_3)\|_{L^2_h}+\|\nabla_hu(.,x_3)\|_{L^2_h}).$$\end{proof}

\noindent Let us now prove Lemma \ref{u,v,1,s}. Bony's decomposition implies
that$$\Delta_q^{\vee}(uv)=T_{1,q}+T_{2,q}+T_{3,q},$$where
\begin{eqnarray*}
T_{1,q}&=&\sum_{|q'-q|\leq
4}\Delta_q^\vee(S^\vee_{q'-1}u.\Delta_{q'}^\vee v),\\T_{2,q}&=&
\sum_{|q'-q|\leq 4}\Delta_q^\vee(S^\vee_{q'-1} v
.\Delta_{q'}^\vee u),\\T_{3,q}&=&\sum_{{\tiny{\begin{array}{c}i\in\{0,\pm 1\}\\  q'>q-4\\
\end{array}}}}\Delta_q^\vee(\Delta_{q'}^\vee u.\Delta_{q'-i}^\vee
v).
\end{eqnarray*}
Proposition \ref{vic} leads to
$$\|S^\vee_{q'-1}u(.,x_3).\Delta_{q'}^\vee v(.,x_3)\|_{L^2_h}\leq C (\|S^\vee_{q'-1}u(.,x_3)
\|_{L^2_h}+\|\nabla_h S^\vee_{q'-1}u(.,x_3)
\|_{L^2_h})\|\Delta_{q'}^\vee v(.,x_3)\|_{L^2_h}.
$$Taking the
 $L^2_{v}$ norm
implies that$$\|S^\vee_{q'-1}u.\Delta_{q'}^\vee v \|_{L^2}\leq C
(\|S^\vee_{q'-1}u\|_{L^{\infty}_v(L^2_h)}+\|\nabla_hS^\vee_{q'-1}u\|_{L^{\infty}_v(L^2_h)})\|\Delta_{q'}^\vee
v\|_{L^2}.$$By Proposition \ref{DP}, we get $$\|T_{1,q}\|_{L^2}\leq
C c_q 2^{-qs} \|v\|_{H^{0,s}}(\|u\|_{H^{0,s}}+\|\nabla_h
u\|_{H^{0,s}}).$$ For the term $T_{2,q}$, Proposition \ref{vic}
leads to
$$\|S^\vee_{q'-1}v(.,x_3).\Delta_{q'}^\vee u(.,x_3)\|_{L^2_h}\leq C \|S^\vee_{q'-1}v(.,x_3)
\|_{L^2_h}(\|\Delta_{q'}^\vee u(.,x_3)\|_{L^2_h}+
\|\nabla_h\Delta_{q'}^\vee u(.,x_3)\|_{L^2_h}).
$$
Taking the $L^2_{v}$ norm implies
that$$\|S^\vee_{q'-1}v.\Delta_{q'}^\vee u \|_{L^2}\leq C
\|S^\vee_{q'-1}v\|_{L^{\infty}_v(L^2_h)}(\|\Delta_{q'}^\vee
u\|_{L^2}+ \|\nabla_h\Delta_{q'}^\vee u\|_{L^2}).$$ Thanks to
Proposition \ref{DP}, we get$$\|T_{2,q}\|_{L^2}\leq C c_q 2^{-qs}
\|v\|_{H^{0,s}}(\|u\|_{H^{0,s}}+\|\nabla_h u\|_{H^{0,s}}).$$ The
term $T_{3,q}$ is estimated as follows. Thanks to the Bernstein
lemma, we get
$$\|T_{3,q}\|_{L^2}\leq C \sum_{{\tiny{\begin{array}{c}i\in\{0,\pm 1\}\\  q'>q-4\\
\end{array}}}}2^{\frac{q}{2}}\|\Delta_{q'}^\vee u.\Delta_{q'-i}^\vee
v\|_{L^1_v(L^2_h)}.$$ Proposition \ref{vic} implies that
$$\|\Delta_{q'}^\vee u(.,x_3).\Delta_{q'-i}^\vee v(.,x_3)\|_{L^2_h}\leq C \|\Delta_{q'-i}^\vee v(.,x_3)\|_{L^2_h}
(\|\Delta_{q'}^\vee u(.,x_3)\|_{L^2_h}+\|\Delta_{q'}^\vee \nabla_h
u(.,x_3)\|_{L^2_h}).$$Taking the $L^1_v$ norm yields

$$\|\Delta_{q'}^\vee u.\Delta_{q'-i}^\vee
v\|_{L^1_v(L^2_h)}\leq C \|\Delta_{q'-i}^\vee v\|_{L^2}
(\|\Delta_{q'}^\vee u\|_{L^2}+\|\Delta_{q'}^\vee \nabla_h
u\|_{L^2}).$$As $s>\frac{1}{2}$, Proposition \ref{DP} implies that
$$\|T_{2,q}\|_{L^2}\leq C c_q 2^{-qs}
\|v\|_{H^{0,s}}(\|u\|_{H^{0,s}}+\|\nabla_h u\|_{H^{0,s}}).$$

Finally, we sum all the above estimates to conclude the proof of
Lemma \ref{u,v,1,s}.
\end{proof}
%@@@@@@@@@@@@@@@@@@@@@@@@@@@@@@@@@@@@%@@@@@@@@@@@@@@@@@@@@@@@@@@@@@@@@@@@@%@@@@@@@@@@@@@@@

\bigskip

\noindent{\bf Acknowledgements}: {\rm The author would like to thank Mohamed Majdoub for suggesting
her this problem and for his kind advice}.

\bigskip

%@@@@@@@@@@@@@@@@@@@@@@@@@@@@@@@@@@@@%@@@@@@@@@@@@@@@@@@@@@@@@@@@@@@@@@@@@%@@@@@@@@@@@@@@@

\end{document}